\documentclass[12pt]{article}
\usepackage{enumerate}
\usepackage{amssymb,a4wide,latexsym,makeidx,epsfig,fleqn}
\usepackage{amsmath}
\usepackage{amsthm}
\usepackage{graphicx}
\usepackage{float}
\usepackage{color}
\usepackage{indentfirst}
\usepackage{pdflscape}
\usepackage{siunitx}
\usepackage{diagbox}
\usepackage{rotating}
\usepackage{booktabs}
\usepackage{titling}
\usepackage[numbers,sort&compress]{natbib}
\usepackage[colorlinks, linkcolor=blue, anchorcolor=blue, citecolor=blue]{hyperref}

\allowdisplaybreaks[4]
\newtheorem{theorem}{Theorem}[section]

\newtheorem{lemma}[theorem]{Lemma}

\newtheorem{claim}{Claim}
\newtheorem*{claim*}{Claim}

\begin{document}
	
	\textwidth 150mm \textheight 225mm
	\topmargin -1.4 cm  \parskip 4pt
	
	\title{Spectral conditions for the existence of $\{K_2,\{C_k\}\}$-factor in $t$-tough graphs\thanks{Supported by the National Natural Science Foundation of China (No. 12271439)}}
	
	\author{{Zhenhao Zhang$^{a,b}$, Ligong Wang$^{a,b,}$\thanks{Corresponding author.}}\\
		{\small $^a$School of Mathematics and Statistics, Northwestern
			Polytechnical University,}\\ {\small  Xi'an, Shaanxi 710129,
			P.R. China.}\\
		{\small $^b$Xi'an-Budapest Joint Research Center for Combinatorics, Northwestern
			Polytechnical University,}\\
		{\small Xi'an, Shaanxi 710129,
			P.R. China. }\\
		{\small E-mails: zhangzhenhao@mail.nwpu.edu.cn, lgwangmath@163.com} }
	\date{}
	\maketitle	
	\vspace{-1cm}
	\begin{center}
		\begin{minipage}{120mm}
			
			\begin{center}
				{\small {\bf Abstract}}
			\end{center}
			{\small 
				For $t>0$, a graph $G$ is called $t$-tough if $|S|\geq t\,c(G-S)$ for every $S\subseteq V(G)$ with $c(G-S)\geq2$, where $c(G-S)$ denotes the number of components of $G-S$. Toughness is closely related to the existence of Hamiltonian cycles and factors in graphs. Chv\'atal [Discrete Math. 5 (1973) 215--228] conjectured that there exists a constant $t_0>0$ such that every $t_0$-tough graph of order at least three is Hamiltonian. Bauer, Broersma and Veldman [Discrete Appl. Math. 99 (2000) 317--321] constructed $(9/4-\varepsilon)$-tough graphs without a Hamiltonian path for arbitrarily small $\varepsilon>0$. For the weaker requirement of a $2$-factor, Enomoto et al. [J. Graph Theory 9 (1985) 87--95] proved that every $2$-tough graph of order at least three has a $2$-factor and that the toughness threshold $2$ is best possible. A graph $G$ is said to have a $\{K_2,\{C_k:k\geq3\}\}$-factor if and only if $G$ contains a spanning subgraph consisting of components that are either cycles or single edges. A $\{K_2,\{C_k:k\geq3\}\}$-factor can be viewed as a combination of a $1$-factor and a $2$-factor, motivating the study of such factors in tough graphs. In this paper, we establish sharp sufficient conditions in terms of the adjacency spectral radius, the distance spectral radius, and the distance signless Laplacian spectral radius for a $t$-tough graph to have a $\{K_2,\{C_k:k\geq3\}\}$-factor.

				 }
		\end{minipage}
	\end{center}
	
	%%%%%%%%%%%%%%%%%%%%%%%%%%%%%%%%%%%%%%%%%%%%%%%%%%%%%%%%%%%%%%%%%%%%%%%%%%%%%%
	\section{Introduction }
	\label{sec:ch6-introduction}
	\noindent

	Throughout this paper, all graphs are assumed to be simple, connected and undirected. Let $G$ be a graph with vertex set $V(G)$ and edge set $E(G)$. The order of $G$ is the cardinality $|V(G)|$. For each vertex $v \in V(G)$, its neighbourhood is denoted by $N_G(v)$, and its degree is given by $d_G(v)=|N_G(v)|$. For a subset  $S \subseteq V(G)$, define $N_G(S)=\bigcup_{v\in S}N_G(v)$. The minimum degree of $G$ is denoted by $\delta(G)$. A vertex of degree zero is called isolated. We write $I(G)$ for the set of isolated vertices of $G$ and $i(G)=|I(G)|$. Let $c(G)$ denote the number of components of a graph $G$.
	Let $E_G(v)=\{e \in E(G)\mid \text{$e$ is incident with  $v$ in $G$} \}$.
	Let $G$ and $H$ be two graphs. If $V(H)\subseteq V(G)$ and $E(H)\subseteq E(G)$, then $H$ is called a subgraph of $G$, denoted by $H \subseteq G$. In addition, if $V(H)=V(G)$ and $H\subseteq G$, then $H$ is a spanning subgraph of $G$, denoted by $H\preceq G$.
	Given a subset $V_1 \subseteq V(G)$, the notation $G-V_1$ refers to the graph obtained by removing the vertices in $V_1$ together with all edges incident with them. For two graphs $G_1$ and $G_2$, their union is denoted by $G_1+G_2$ and has vertex set $V(G_1)\cup V(G_2)$ and edge set $E(G_1)\cup E(G_2)$. When $G_1$ and $G_2$ are vertex-disjoint, their join, denoted by $G_1\vee G_2$, is formed from their union by connecting every vertex of $G_1$ to every vertex of $G_2$. The complement of $G$, denoted by $\overline{G}$, is a graph on the same vertex set $V(G)$ such that two distinct vertices $u$ and $v$ are adjacent in $\overline{G}$ if and only if they are not adjacent in $G$.
	Let $K_n$ and $C_n$ denote the complete graph and the cycle of order $n$, respectively.
	
	For a graph $G$ of order $n$, let $A(G)=(a_{i,j})_{n\times n}$ be its adjacency matrix, where $a_{i,j}=1$ when $v_i$ and $v_j$ are adjacent and $a_{i,j}=0$ otherwise.
	For two vertices $v_i,v_j\in V(G)$, let $d_G(v_i,v_j)$, or simply $d_{i,j}$, denote the distance between $v_i$ and $v_j$ in $G$, that is, the length of a shortest path joining them. The distance matrix of $G$ is the real symmetric matrix $D(G)=(d_{i,j})$. For a vertex $v_i\in V(G)$, its transmission in $G$, denoted by $\text{Tr}_G(v_i)$, is defined as the sum of the distances from $v_i$ to all other vertices of $G$. Let $\text{Tr}(G)=\text{diag}\big(\text{Tr}_G(v_1),\text{Tr}_G(v_2),\ldots,\text{Tr}_G(v_n)\big)$ be the diagonal matrix whose diagonal entries are the transmissions of the vertices of $G$. The distance signless Laplacian matrix of $G$ is then defined by $Q(G)=\text{Tr}(G)+D(G)$. The largest eigenvalues of $A(G)$, $D(G)$ and $Q(G)$ are called the adjacency spectral radius, the distance spectral radius and the distance signless Laplacian spectral radius of $G$, respectively, and are denoted by $\rho(G)$,
	$\mu(G)$ and $\eta(G)$.
	
	Given subgraphs $H_1, H_2,\ldots,H_k$ of $G$,  a $\{H_1, H_2,\ldots,H_k\}$-factor of   $G$ is a spanning subgraph $F$ in which each connected component is isomorphic to one of $H_1, H_2,\ldots,H_k$. A graph $G$ is said to have a $\{K_2,\{C_{k}:k\geq 3\}\}$-factor if and only if $G$ contains a spanning subgraph consisting of components that are either cycles or single edges. For brevity, we refer to this as a $\{K_2,\{C_k\}\}$-factor. Tutte \cite{Tut} established a necessary and sufficient condition for the existence of a $\{K_2,\{C_k\}\}$-factor in a graph from a structural perspective. Motivated by this result, Lou et al. \cite{Lou} gave a tight sufficient condition in terms of the spectral radius for the existence of a $\{K_2,\{C_k\}\}$-factor in a graph with minimum degree $\delta$. Xu et al. \cite{Xu} extended the above result to the distance spectral radius of a graph. For more results on factors, the reader is referred to the surveys \cite{Plummer} and \cite{Fan1}.
	
	Let $t$ be a positive real number. A graph $G$ is called $t$-tough if $|S| \geq t\,c(G-S)$
	for every subset $S\subseteq V(G)$ with $c(G-S)\geq 2$. The toughness of a non-complete graph $G$, denoted by $t(G)$, is defined as
	\[
	t(G)=\min\left\{
	\frac{|S|}{c(G-S)}:
	S\subseteq V(G),\ c(G-S)\geq 2
	\right\}.
	\]
	For the complete graph $K_n$, we set $t(K_n)=+\infty$.

	The toughness of a graph is closely related to its spectral radius.	Gu and Haemers \cite{Gu} obtained two lower bounds for the toughness of a graph in terms of its Laplacian eigenvalues. Fan, Lin and Lu \cite{Fan} subsequently established adjacency spectral conditions for a graph with prescribed minimum degree to be $1$-tough, and for a graph to be $r$-tough when $r\geq1$ is an integer. They also characterized the corresponding extremal graphs. For a variant of toughness introduced by Enomoto, Chen, Fan and Lin \cite{Chen} gave further adjacency spectral conditions for a graph with minimum degree $\delta$ to be $r$-tough and $t$-tough, respectively, where $r\geq 2$ and $1/t$ are integers. They extended their investigation to balanced bipartite graphs. More recently, Zhai et al. \cite{Zhai} determined the maximum adjacency spectral radii for several classes of tough graphs, including balanced bipartite classes, and identified the corresponding extremal graphs. For the distance spectral radius, Lou, Liu and Shu \cite{Lou1} established a sufficient condition  in terms of the distance spectral radius for a graph with given minimum degree to be $1$-tough. They also obtained distance spectral conditions for $t$-toughness when $t$ or $1/t$ is a positive integer. Toughness is also closely related to the existence of factors. Chv\'atal \cite{Chvatal} introduced  the concept of toughness and conjectured that there exists a constant $t_0>0$ such that every $t_0$-tough graph of order at least three is Hamiltonian. Bauer, Broersma and Veldman \cite{Bauer} constructed
	$(9/4-\varepsilon)$-tough graphs without a Hamiltonian path for arbitrarily small $\varepsilon>0$, implying that such a constant, if it exists, must be at least $9/4$. Since a Hamiltonian cycle is a connected $2$-factor, a natural weaker question is whether
	sufficiently large toughness guarantees the existence of a $2$-factor. Enomoto et al. \cite{Enomoto} proved that every $2$-tough graph of order at least three has a $2$-factor. Moreover, the threshold $2$ is sharp: for every $0<\varepsilon<2$, there exists a $(2-\varepsilon)$-tough graph with no $2$-factor. A $\{K_2,\{C_k:k\geq3\}\}$-factor can be viewed as a combination of a $1$-factor and a $2$-factor. Motivated by investigations on such factors in tough graphs, we derive sufficient conditions involving the adjacency spectral radius, the distance spectral radius, and the distance signless Laplacian spectral radius for a $t$-tough graph to have a $\{K_2,\{C_k\}\}$-factor.

	The remainder of this paper is organized as follows. Section 2 introduces some useful lemmas that will be used in the subsequent proofs. In Section 3, we establish two adjacency spectral radius conditions for a $t$-tough graph to have a $\{K_2,\{C_k\}\}$-factor. In Section 4, we establish two distance spectral radius conditions for a $t$-tough graph to have a $\{K_2,\{C_k\}\}$-factor. Finally, in Section 5, we establish two distance signless Laplacian spectral radius conditions for a $t$-tough graph to have a $\{K_2,\{C_k\}\}$-factor.

	\section{Preliminaries}\label{sec:Preliminaries}
	
	In this section, we give some useful lemmas which will be used later.

	\begin{lemma}\label{lem0}
	(\cite{Tut1, Lou})
	Let $G$ be a graph and let $k \geq 3$ be an integer. $G$ has a $\{K_2,\{C_k\}\}$-factor if and only if $i(G-S)\leq |S|$ for every $ S\subseteq V(G)$.
	
    \end{lemma}

    For a square matrix $M$, let $\rho(M)$ denote its spectral radius. Given two $n\times n$ matrices $A=(a_{i,j})$ and $B=(b_{i,j})$, we write $A\geq B$ whenever $a_{i,j}\geq b_{i,j}$ holds for all $i,j\in\{1,2,\dots,n\}$. The matrix $A=(a_{i,j})$ is said to be nonnegative if $a_{i,j}\geq 0$ for all $i,j\in\{1,2,\dots,n\}$.
	
	\begin{lemma}\label{lem5}	\emph{(\cite{Ber})}
		Let $A=(a_{i,j})$ and $B=(b_{i,j})$ be two nonnegative square matrices of order $n$. If $A\geq B$, then $\rho(A)\geq\rho(B)$. Moreover, if $A$ is irreducible, $A\geq B$, and $A\neq B$, then $\rho(A)>\rho(B)$.
	\end{lemma}
	
	Let $M$ be a square matrix of order $n$, and let $\mathcal{N}=\{1,2,\dots,n\}$. Suppose
	$\Pi:\mathcal{N}=\mathcal{N}_1\cup \mathcal{N}_2\cup \cdots \cup \mathcal{N}_k$
	is a partition of $\mathcal{N}$. For $i,j\in\{1,2,\dots,k\}$, denote by $M_{i,j}$ the submatrix of $M$ formed by the rows indexed by $\mathcal{N}_i$ and the columns indexed by $\mathcal{N}_j$. Then $M$ is similar to the block matrix $M'$ given by
	\[
	M'=\begin{pmatrix}
		M_{1,1}&M_{1,2}&\cdots&M_{1,k}\\
		M_{2,1}&M_{2,2}&\cdots&M_{2,k}\\
		\vdots&\vdots&\ddots&\vdots\\
		M_{k,1}&M_{k,2}&\cdots&M_{k,k}
	\end{pmatrix}.
	\]
	
	With respect to $\Pi$, the quotient matrix of $M$, denoted by $B_{\Pi}=(b_{i,j})$, is the $k\times k$ matrix whose entry $b_{i,j}$ is the average row sum over the block $M_{i,j}$. The partition $\Pi$ is called equitable if every block $M_{i,j}$ has constant row sum; in this case, the quotient matrix $B_{\Pi}$ is also said to be equitable.
	
	\begin{lemma}\label{lem6}	\emph{(\cite{Bro, God})}
		Let $M$ be a real symmetric matrix and $B_{\Pi}$ an equitable quotient matrix of $M$. Then every eigenvalue of $B_{\Pi}$ is an eigenvalue of $M$. Moreover, if $M$ is nonnegative and irreducible, then the largest eigenvalue of $B_{\Pi}$ coincides with the largest eigenvalue of $M$.
	\end{lemma}
	
	\begin{lemma}\label{lem11} \emph{(\cite{Ber})}
		Let $M$ be an irreducible nonnegative matrix of order $n$ and let $x=(x_1,\ldots,x_n)^{\mathsf T}>0$. Then
		\[\min_{1\leq i\leq n}\frac{(Mx)_i}{x_i}\leq\rho(M)\leq\max_{1\leq i\leq n}\frac{(Mx)_i}{x_i}.\]
		Equality in either inequality holds if and only if $Mx=\rho(M)x$. 
		Moreover, for any $\lambda\in\mathbb{R}$, if $Mx\leq\lambda x$ and $Mx\neq\lambda x$, then $\rho(M)<\lambda$; if $Mx\geq\lambda x$ and $Mx\neq\lambda x$, then $\rho(M)>\lambda$. All vector inequalities are understood componentwise.
	\end{lemma}
	
	\begin{lemma}\label{lem7}\emph{(\cite{Lou}, Theorem~1.5, Case~2)}
		Let $G$ be a graph of order $n \geq 5a + 6$ where $a\geq 1$ is an integer. Denote $G_s=K_s\vee (K_{n-2s-1}+(s+1){K_1})$ with $a\leq s\leq \lfloor \frac{n-3}{2}\rfloor$. Then $ \rho(G_s)\leq \rho(G_{a})$ with equality if and only if $s=a$.
	\end{lemma}
	
	\begin{lemma}\label{lem8}\emph{(\cite{Xu}, Theorem~4.1, Case~2)}
		Let $G$ be a graph of order $n \geq 10a + 11$ where $a\geq 1$ is an integer. Denote $G_s=K_s\vee (K_{n-2s-1}+(s+1){K_1})$ with $a\leq s\leq \lfloor \frac{n-3}{2}\rfloor$. Then $ \mu(G_s)\geq \mu(G_{a})$ with equality if and only if $s=a$.
	\end{lemma}
	
	\begin{lemma}\label{lem9}		
		Let $G$ be a graph of order $n \geq 14a + 4$ where $a\geq 1$ is an integer. Denote $G_s=K_s\vee (K_{n-2s-1}+(s+1){K_1})$ with $a\leq s\leq \lfloor \frac{n-3}{2}\rfloor$. Then $ \eta(G_s)\geq \eta(G_{a})$ with equality if and only if $s=a$.
	\end{lemma}
	
	\begin{proof}
   The equality is immediate when $s=a$. Suppose that $s>a$.
	With the partition $V(G_s) = V(K_{n-2s-1}) \cup V({K_s})\cup V((s+1)K_1)$, the corresponding quotient matrix of $Q(G_s)$ is
	\[
	M_s=
	\begin{pmatrix}
		2n-s-2&s&2(s+1)\\
		n-2s-1&n+s-2&s+1\\
		2(n-2s-1)&s&2n+s-2
	\end{pmatrix}.
	\]
	
	The characteristic polynomial of $M_s$ is
	\begin{align*}
		f_s(x)=&|xI_3-M_s|\\
		=&x^3+(6-5n-s)x^2+(8n^2-ns-24n+8s^2+8s+16)x\\
		-&4n^3+2n^2s+20n^2-8ns^2-14ns-32n-2s^3+14s^2+20s+16.
	\end{align*}
	
	The characteristic polynomial of $M_a$ is
	\begin{align*}
		f_a(x)=&|xI_3-M_a|\\
		=&x^3+(6-5n-a)x^2+(8n^2-na-24n+8a^2+8a+16)x\\
		-&4n^3+2n^2a+20n^2-8na^2-14na-32n-2a^3+14a^2+20a+16.
	\end{align*}
	
	Let $\theta_{s}$ and $\theta_{a}$ be the largest roots of $f_s(x)$ and $f_a(x)$, respectively. Hence $\theta_{s}=\eta(G_s)$ and $\theta_{a}=\eta(G_a)$.
	Let $g(x)=f_s(x)-f_a(x)$. Then  $g(x)=(a-s)[x^2+(n-8s-8a-8)x+2s^2+2as+2a^2+8n(s+a)-14(s+a)-2n^2+14n-20]$.
	Notice that $\eta(G_a)\geq \eta(K_n)=2n-2$. We claim that
	$g(x)<0$ for $x \in [ 2n-2, +\infty )$.
	
	Let $g(x)=(a-s) \widehat{g}(x)$, where
    \begin{align*}
		\widehat g(x)
	&=x^2+(n-8s-8a-8)x+2s^2+2as+2a^2\\
	&\quad+8n(s+a)-14(s+a)-2n^2+14n-20.
	\end{align*}
    Since $n\geq 14a+4$ and $a< s\leq \lfloor\frac{n-3}{2}\rfloor \leq \frac{n-3}{2}$,
    \[
    2n-2-(-\frac{n-8s-8a-8}{2})=\frac{5n}{2}-4s-4a-6 \geq \frac{n}{2}-4a\geq 0.
    \]
    Hence $\widehat{g}(x)$ is increasing on the interval $ [ 2n-2, +\infty )$. Let $h(s)=\widehat{g}(2n-2)$. Then
    \begin{align*}
    	h(s)
    	={}&2s^2+(2a-8n+2)s\\
    	&+4n^2-8an-12n+2a^2+2a.
    \end{align*}
    As $s \leq \lfloor \frac{n-3}{2}\rfloor \leq \frac{n-3}{2}$, we have
    \[
    h'(s)=4s+2a-8n+2\leq2a-6n<0.
    \]
    Consequently,
    \begin{align*}
    	h(s)
    	&\geq h\left(\frac{n-3}{2}\right)\\
    	&=n\left(\frac n2-7a-2\right)
    	+2a^2-a+\frac32>0.
    \end{align*}

    Therefore $g(x)=(a-s)\widehat g(x)\leq (a-s)h(s)<0$ for $x\geq2n-2$. $0>g(\theta_a)=f_s(\theta_a)-f_a(\theta_a)=f_s(\theta_a)$. Hence  $\eta(G_s)> \eta(G_a)$ when $s>a$. When $s=a$, we have $G_s=G_a$, and thus $\eta(G_s)=\eta(G_a)$. Therefore $\eta(G_s)\geq\eta(G_a)$, with equality if and only if $s=a$.
    This completes the proof.
	\end{proof}

	\begin{lemma}\label{lem10}		
	Let $0<t<1$.
		
	(1) Let $n,s,p$ be positive integers with $n\geq s+p+2$, and let $H_1=K_s\vee(K_{n-s-p}+pK_1)$. If $t(H_1)\geq t$ and $H_1$ has no $\{K_2,\{C_k\}\}$-factor, then $t\leq\frac{2\lceil\frac{n}{2}\rceil-4}{2\lceil\frac{n}{2}\rceil}$.
		
	(2) Let $n,p$ be positive integers with $1\leq p\leq n-1$, and let $H_2=K_{n-p}\vee pK_1$. If $t(H_2)\geq t$ and $H_2$ has no $\{K_2,\{C_k\}\}$-factor, then $t\leq\frac{2\lceil\frac{n}{2}\rceil-2}{2\lfloor\frac{n}{2}\rfloor+2}$.
    \end{lemma}
    
    \begin{proof}
    	(1) By Lemma \ref{lem0}, since $H_1$ does not have a $\{K_2,\{C_k\}\}$-factor, there exists $S_1 \subseteq V(H_1)$ such that  $i(H_1-S_1)-|S_1|>0$. We claim that $i(H_1-V(K_s))-s>0$, that is, $p>s$. Assume, to the contrary, that $p\leq s$. Suppose  $V(K_s)\not \subseteq S_1$. Then a vertex of $K_s$ remains in $H_1-S_1$ and is adjacent to every other remaining vertex. Consequently, $H_1-S_1$ is either connected or consists of a single vertex. In both cases, $i(H_1-S_1)\leq 1 \leq |S_1|$.  It contradicts the choice of $S_1$. Hence, $V(K_s)\subseteq S_1$. When $S_1=V(K_s)$, then $i(H_1-S_1)-|S_1|= p-s \leq 0$. When $V(K_s)\subset S_1$, then $|S_1|\geq s+1$ and $i(H_1-S_1)\leq p+1$. Hence $i(H_1-S_1)-|S_1|\leq p+1-(s+1)=p-s\leq 0$. In either case, we obtain a contradiction. Therefore, $p>s$. 
    	
    	Since $t(H_1)\geq t$, we have $\frac{|S_1|}{c(H_1-S_1)}\geq t$ for every vertex cut    $S_1\subseteq V(H_1)$ with $c(H_1-S_1)>1$. In particular, taking $S_1=V(K_s)$, we obtain $s\geq t(p+1)$. Thus,
    	\begin{align*}
    		\begin{cases}
    			p\geq s+1 \\
    			s\geq t(p+1) \\
    			n-s-p\geq 2
    		\end{cases} 
    		\quad \Rightarrow \quad
    		\begin{cases}
    			s\geq t(s+2)\\
    			n-2s-1\geq 2
    		\end{cases}
    		\quad \Rightarrow \quad   
    		\lceil \frac{2t}{1-t}\rceil \leq s \leq \lfloor \frac{n-3}{2}\rfloor.
    	\end{align*}
    	
    	It means that $\lceil \frac{2t}{1-t}\rceil  \leq \lfloor \frac{n-3}{2}\rfloor$.
    	Let $B_1=\lfloor \frac{n-3}{2}\rfloor$. Then $B_1$ is an integer. An integer $B_1\geq \lceil \frac{2t}{1-t}\rceil$ if and only if $B_1\geq   \frac{2t}{1-t}$.

    	\begin{align*}
    		\begin{cases}
    			\lfloor \frac{n-3}{2}\rfloor= \frac{n-3}{2} \geq   \frac{2t}{1-t}, &\text{when $n$ is odd,}\\
    			\lfloor \frac{n-3}{2}\rfloor=\frac{n}{2}-2\geq   \frac{2t}{1-t}, &\text{when $n$ is even.}
    		\end{cases} 
    		\quad \Rightarrow \quad
    		\begin{cases}
    			t\leq  \frac{n-3}{n+1},&\text{when $n$ is odd},\\
    			t\leq \frac{n-4}{n},&\text{when $n$ is even}.
    		\end{cases}
    	\end{align*}
    	
    	Combining the above situations, we obtain $t \leq \frac{2\lceil \frac{n}{2}\rceil-4}{2\lceil  \frac{n}{2}\rceil}$.
    	
    	(2) By Lemma \ref{lem0}, since $H_2$ does not have a $\{K_2,\{C_k\}\}$-factor, there exists $S_2 \subseteq V(H_2)$ such that  $i(H_2-S_2)-|S_2|>0$. We claim that $i(H_2-V(K_{n-p}))-(n-p)>0$, that is, $p\geq n-p+1$. Assume, to the contrary, that $p\leq n-p$. Suppose $V(K_{n-p})\not \subseteq S_2$. Then a vertex of $K_{n-p}$ remains in $H_2-S_2$ and is adjacent to every other remaining vertex. Consequently, $H_2-S_2$ is either connected or consists of a single vertex. In both cases, $i(H_2-S_2)\leq 1 \leq |S_2|$.  It contradicts the choice of $S_2$. Hence, $V(K_{n-p})\subseteq S_2$.  $i(H_2-S_2)-|S_2|\leq p-(n-p)\leq 0$. This leads to a contradiction. Therefore, $p\geq n-p+1$. 
    	
    	Taking $S_2=V(K_{n-p})$, $\frac{|S_2|}{c(H_2-S_2)}=\frac{n-p}{p}\geq t(H_2)\geq t$. Thus,

    		\begin{align*}
    		\begin{cases}
    			p\geq n-p+1 \\
    			n-p\geq tp 
    		\end{cases} 
    		\quad \Rightarrow \quad
    		\begin{cases}
    			p\geq  \lceil \frac{n+1}{2}\rceil,\\
    			p\leq \lfloor  \frac{n}{1+t}\rfloor.
    		\end{cases}
    	\end{align*}
    	
    	It means that $ \lceil \frac{n+1}{2}\rceil \leq \lfloor  \frac{n}{1+t}\rfloor$. Let $B_2=\lceil \frac{n+1}{2}\rceil$. Then $B_2$ is an integer. An integer $B_2\leq \lfloor  \frac{n}{1+t}\rfloor $ if and only if $B_2\leq   \frac{n}{1+t}$.
    	\begin{align*}
    		\begin{cases}
    			\lceil \frac{n+1}{2}\rceil= \frac{n+1}{2} \leq   \frac{n}{1+t}, &\text{when $n$ is odd,}\\
    			\lceil \frac{n+1}{2}\rceil=\frac{n}{2}+1\leq   \frac{n}{1+t}, &\text{when $n$ is even.}
    		\end{cases} 
    		\quad \Rightarrow \quad
    		\begin{cases}
    			t\leq  \frac{n-1}{n+1},&\text{when $n$ is odd},\\
    			t\leq \frac{n-2}{n+2},&\text{when $n$ is even}.
    		\end{cases}
    	\end{align*}
    	
    	Combining the above situations, we obtain $t \leq \frac{2\lceil \frac{n}{2}\rceil-2}{2\lfloor  \frac{n}{2}\rfloor+2}$. 
    \end{proof}

	\section{Two adjacency spectral conditions}
	
	In this section, we establish two tight sufficient conditions in terms of the spectral radius for a $t$-tough graph $G$ to have a $\{K_2,\{C_k\}\}$-factor. Theorem  \ref{thm1} gives the spectral radius condition for $G$ to have a $\{K_2,\{C_k\}\}$-factor, which is divided into three cases based on the value of $t$. Furthermore, for graphs of larger order $n \ge \frac{12}{1-t}$, Theorem \ref{thm2} provides a more concise spectral radius condition.
	
	If $t\geq 1$, then for a $t$-tough graph $G$ and every vertex cut $S\subseteq V(G)$,  $\frac{|S|}{c(G-S)}\geq t$. Hence $|S|\geq tc(G-S)\geq ti(G-S)\geq i(G-S)$, and $G$ has a $\{K_2,\{C_k\}\}$-factor. Consequently we may assume $0<t<1$.

	\begin{theorem}\label{thm1}
		Let $0<t<1$ and let $G$ be a $t$-tough graph of order $n$. Set $G_s=K_s\vee (K_{n-2s-1}+(s+1)K_1)$ with  $\lceil \frac{2t}{1-t}\rceil \leq s \leq \lfloor \frac{n-3}{2}\rfloor$. If one of the following three conditions holds: 
		
		\enspace	(1) $\frac{2\lceil \frac{n}{2}\rceil-2}{2\lfloor  \frac{n}{2}\rfloor+2}< t<1$;
		
		\enspace	(2) $\frac{2\lceil \frac{n}{2}\rceil-4}{2\lceil \frac{n}{2}\rceil}< t\leq \frac{2\lceil \frac{n}{2}\rceil-2}{2\lfloor\frac{n}{2}\rfloor+2}$, $\rho(G)\geq \rho(K_{\lfloor  \frac{n-1}{2}\rfloor}\vee \lceil \frac{n+1}{2}\rceil K_1)$, and $G\neq K_{\lfloor  \frac{n-1}{2}\rfloor}\vee \lceil \frac{n+1}{2}\rceil K_1$;
		
		\enspace	(3) $t\leq \frac{2\lceil \frac{n}{2}\rceil-4}{2\lceil  \frac{n}{2}\rceil}$, $\rho(G)\geq \mbox{max}\{\rho(K_{\lfloor  \frac{n-1}{2}\rfloor}\vee \lceil \frac{n+1}{2}\rceil K_1), \rho(G_s)\}$, and $G$ is isomorphic to neither $K_{\lfloor  \frac{n-1}{2}\rfloor}\vee \lceil \frac{n+1}{2}\rceil K_1$ nor any $G_s$ with $\left\lceil\frac{2t}{1-t}\right\rceil \leq s\leq \left\lfloor\frac{n-3}{2}\right\rfloor$;
				
		then $G$ has a $\{K_2,\{C_k\}\}$-factor.
	\end{theorem}

    \begin{proof}

	Assume to the contrary that $G$ does not have a $\{K_2,\{C_k\}\}$-factor. By Lemma \ref{lem0}, there exists $S \subseteq V(G)$ such that  $i(G-S)-|S|>0$. Let $|S|=s$, $i(G-S)=p$ and $c(G-S)=p+q$. Then $p>s$. We now proceed to discuss the following three cases. 
	
	{\textbf{Case 1.}} $\frac{2\lceil \frac{n}{2}\rceil-2}{2\lfloor  \frac{n}{2}\rfloor+2}< t<1$.
	
	If there are nontrivial components in $G-S$, we can add all possible edges between the nontrivial components of $G-S$ so that they form a clique. Then $G\preceq K_s\vee (K_{n-s-p}+p{K_1})$. Set $H_1=K_s\vee (K_{n-s-p}+p{K_1})$. Then $t(H_1)\geq t(G)\geq t$ and $H_1$ does not have a $\{K_2,\{C_k\}\}$-factor. By Lemma \ref{lem10}(1), $t\leq \frac{2\lceil \frac{n}{2}\rceil-4}{2\lceil  \frac{n}{2}\rceil}$. This implies that $t\leq \frac{2\lceil \frac{n}{2}\rceil-4}{2\lceil  \frac{n}{2}\rceil}< \frac{2\lceil \frac{n}{2}\rceil-2}{2\lfloor  \frac{n}{2}\rfloor+2}< t$, which is a contradiction.
	
	If there are no nontrivial components in $G-S$, then $G \preceq K_{n-p}\vee pK_1$. Set $H_2=K_{n-p}\vee pK_1$. Since $t(H_2)\geq t(G)\geq t$ and $H_2$ does not have a $\{K_2,\{C_k\}\}$-factor, by Lemma \ref{lem10}(2), $t\leq \frac{2\lceil \frac{n}{2}\rceil-2}{2\lfloor\frac{n}{2}\rfloor+2}$. This implies that $t\leq \frac{2\lceil \frac{n}{2}\rceil-2}{2\lfloor\frac{n}{2}\rfloor+2}<t$, which is a contradiction. 
	
	Thus, both possibilities lead to a contradiction. The proof of Theorem \ref{thm1}(1) is complete.

	{\textbf{Case 2.}} $\frac{2\lceil \frac{n}{2}\rceil-4}{2\lceil \frac{n}{2}\rceil}< t\leq \frac{2\lceil \frac{n}{2}\rceil-2}{2\lfloor\frac{n}{2}\rfloor+2}$.
	
	If there are nontrivial components in $G-S$, $G\preceq K_s\vee (K_{n-s-p}+p{K_1})$. Set $H_1=K_s\vee (K_{n-s-p}+p{K_1})$. Then $t(H_1)\geq t(G)\geq t$ and $H_1$ does not have a $\{K_2,\{C_k\}\}$-factor. By Lemma \ref{lem10}(1), $t\leq \frac{2\lceil \frac{n}{2}\rceil-4}{2\lceil  \frac{n}{2}\rceil}$. This implies that $t\leq \frac{2\lceil \frac{n}{2}\rceil-4}{2\lceil  \frac{n}{2}\rceil}<t$, which is a contradiction.
	
	If there are no nontrivial components in $G-S$, then $G \preceq K_{n-p}\vee pK_1$. Since $p>s=n-p$, we have $p\geq \lceil \frac{n+1}{2}\rceil$. By moving $p-\lceil \frac{n+1}{2}\rceil$ vertices from the independent set into the clique and adding the required edges, we see that $K_{n-p}\vee pK_1\preceq K_{\lfloor  \frac{n-1}{2}\rfloor}\vee \lceil \frac{n+1}{2}\rceil K_1$. Therefore, $\rho(G)\leq \rho(K_{n-p}\vee pK_1)\leq \rho(K_{\lfloor \frac{n-1}{2}\rfloor}\vee \lceil \frac{n+1}{2}\rceil K_1)$ with equality throughout if and only if $G=K_{\lfloor  \frac{n-1}{2}\rfloor}\vee \lceil \frac{n+1}{2}\rceil K_1$. Together with the hypothesis $\rho(G)\geq\rho(K_{\lfloor  \frac{n-1}{2}\rfloor}\vee \lceil \frac{n+1}{2}\rceil K_1)$, this forces $G=K_{\lfloor  \frac{n-1}{2}\rfloor}\vee \lceil \frac{n+1}{2}\rceil K_1$, contradicting the exclusion in (2). 
	
    The proof of Theorem \ref{thm1}(2) is complete.
	
	{\textbf{Case 3.}} $t\leq \frac{2\lceil \frac{n}{2}\rceil-4}{2\lceil  \frac{n}{2}\rceil}$.
	
	If there are nontrivial components in $G-S$, $G\preceq K_s\vee (K_{n-s-p}+p{K_1})$ and $n-s-p\geq 2$. Since $G$ is $t$-tough and $c(G-S)=p+q\geq p+1$,
	we have $s\geq t(p+q)\geq t(p+1)\geq t(s+2)$. Together with $p\geq s+1$ and $n-s-p\geq 2$, it gives that $\lceil\frac{2t}{1-t}\rceil\leq s\leq\lfloor\frac{n-3}{2}\rfloor$. By moving $p-s-1$ isolated vertices into the clique $K_{n-s-p}$ and adding the required edges, we obtain $H_1\preceq G_s=K_s\vee (K_{n-2s-1}+(s+1)K_1)$. Hence $\rho(G)\leq \rho(H_1)\leq \rho(G_s)\leq \max\{\rho(K_{\lfloor  \frac{n-1}{2}\rfloor}\vee \lceil \frac{n+1}{2}\rceil K_1), \rho(G_s)\}$. If $\rho(G)<\rho(K_{\lfloor  \frac{n-1}{2}\rfloor}\vee \lceil \frac{n+1}{2}\rceil K_1)$, then $\rho(G)<\max\{\rho(K_{\lfloor  \frac{n-1}{2}\rfloor}\vee \lceil \frac{n+1}{2}\rceil K_1), \rho(G_s)\}$, contradicting the hypothesis. If $\rho(G)=\max\{\rho(K_{\lfloor  \frac{n-1}{2}\rfloor}\vee \lceil \frac{n+1}{2}\rceil K_1), \rho(G_s)\}$, the hypothesis $\rho(G)\geq \max\{\rho(K_{\lfloor  \frac{n-1}{2}\rfloor}\vee \lceil \frac{n+1}{2}\rceil K_1), \rho(G_s)\}$ forces equality throughout. Then $G= G_s$ or $G=K_{\lfloor  \frac{n-1}{2}\rfloor}\vee \lceil \frac{n+1}{2}\rceil K_1$, contrary to the exclusion in (3).
	Thus, both possibilities lead to a contradiction.
	
	If there are no nontrivial components in $G-S$, then $G \preceq K_{n-p}\vee pK_1\preceq K_{\lfloor  \frac{n-1}{2}\rfloor}\vee \lceil \frac{n+1}{2}\rceil K_1$. Hence $\rho(G)\leq \rho(K_{\lfloor  \frac{n-1}{2}\rfloor}\vee \lceil \frac{n+1}{2}\rceil K_1)\leq \max\{\rho(K_{\lfloor  \frac{n-1}{2}\rfloor}\vee \lceil \frac{n+1}{2}\rceil K_1),  \rho(G_s)\}$. If $\rho(G)<\rho(K_{\lfloor  \frac{n-1}{2}\rfloor}\vee \lceil \frac{n+1}{2}\rceil K_1)$, then $\rho(G)<\max\{\rho(K_{\lfloor  \frac{n-1}{2}\rfloor}\vee \lceil \frac{n+1}{2}\rceil K_1), \rho(G_s)\}$, contradicting the hypothesis. If $\rho(G)=\max\{\rho(K_{\lfloor  \frac{n-1}{2}\rfloor}\vee \lceil \frac{n+1}{2}\rceil K_1), \rho(G_s)\}$, the hypothesis $\rho(G)\geq \max\{\rho(K_{\lfloor  \frac{n-1}{2}\rfloor}\vee \lceil \frac{n+1}{2}\rceil K_1), \rho(G_s)\}$ forces equality throughout. Then $G= G_s$ or $G=K_{\lfloor  \frac{n-1}{2}\rfloor}\vee \lceil \frac{n+1}{2}\rceil K_1$, contrary to the exclusion in (3).
	Thus, both possibilities lead to a contradiction.
	
    The proof of Theorem \ref{thm1}(3) is complete.
    \end{proof}
   
    \begin{theorem}\label{thm2}
    Let $G$ be a $t$-tough graph of order $n\geq \frac{12}{1-t}$, where $0<t<1$. If $\rho(G)\geq \rho (K_{\lceil\frac{2t}{1-t}\rceil} \vee (K_{n-2\lceil\frac{2t}{1-t}\rceil-1}+(\lceil\frac{2t}{1-t}\rceil+1) K_1))  $, then $G$ has a $\{K_2,\{C_k\}\}$-factor unless $G=K_{\lceil\frac{2t}{1-t}\rceil} \vee (K_{n-2\lceil\frac{2t}{1-t}\rceil-1}+(\lceil\frac{2t}{1-t}\rceil+1) K_1)$. 
   
    \end{theorem}
   
    \begin{proof}
   	Assume to the contrary that $G$ does not have a $\{K_2,\{C_k\}\}$-factor and is not isomorphic to the exceptional graph stated in the theorem. By Lemma \ref{lem0}, there exists $S \subseteq V(G)$ such that  $i(G-S)-|S|>0$. Let $|S|=s$, $i(G-S)=p$ and $c(G-S)=p+q$. Then $p>s$. We now proceed to discuss the following two cases. 
   	
   	{\textbf{Case 1.}} $q>0$.
   	
   	Suppose that the nontrivial components of $G-S$ are $G^{1}, G^{2}, \ldots,G^{q}$ and $|V(G^1)|=n_1, |V(G^2)|=n_2, \ldots, |V(G^q)|=n_q$ with $n_1 \geq n_2 \geq \cdots \geq n_q \geq 2 $. Then $G$ is a spanning subgraph of $K_s \vee (K_{n_1}+K_{n_2}+\cdots+K_{n_q}+pK_1)$. We can add all possible edges so that the nontrivial components form a clique in $G-S$. Then $G \preceq K_s \vee (K_{n-s-p} +pK_1)$. Moreover, since $p>s$, we can move $p-s-1$ isolated vertices into the clique $K_{n-s-p}$ and add the required edges. Hence $G\preceq K_s \vee (K_{n-s-p} +pK_1)\preceq K_s\vee (K_{n-2s-1}+(s+1)K_1)$. Since $G$ is $t$-tough and $c(G-S)=p+q$,
   	we have $s\geq t(p+q)\geq t(p+1)\geq t(s+2)$. Together with $p\geq s+1$ and $n-s-p\geq 2$, it gives that $\lceil\frac{2t}{1-t}\rceil\leq s\leq\lfloor\frac{n-3}{2}\rfloor$. Denote $G_s=K_s\vee (K_{n-2s-1}+(s+1){K_1})$. Then $\rho(G)\leq \max\{\rho(G_s): \lceil \frac{2t}{1-t}\rceil \leq s \leq \lfloor \frac{n-3}{2}\rfloor \}$ in this case. 
   	By Lemma \ref{lem7}, we have $\rho(G)\leq \rho(G_s)\leq \rho(G_{\lceil \frac{2t}{1-t}\rceil})$ with equality if and only if $G=G_{\lceil \frac{2t}{1-t}\rceil}$ and $s=\lceil \frac{2t}{1-t}\rceil$. 
   	On the other hand, the hypothesis of Theorem \ref{thm2} requires $\rho(G)\geq \rho (G_{\lceil \frac{2t}{1-t}\rceil})$. Consequently, $\rho(G)=\rho(G_{\lceil \frac{2t}{1-t}\rceil})$, and hence $G=G_{\lceil \frac{2t}{1-t}\rceil}$, contradicting the assumption that $G$ is not the exceptional graph.

   	{\textbf{Case 2.}} $q=0$.
   	
   	In this case, we have $G\preceq K_{n-p}\vee p{K_1}$. Hence $\rho(G)\leq \rho(K_{n-p}\vee p{K_1})\leq \rho(K_{\lfloor \frac{n-1}{2}\rfloor}\vee \lceil \frac{n+1}{2}\rceil{K_1})$. We then state the following claim.
   	
   	\begin{claim} \label{claim1}
   		$\rho(K_{\lfloor \frac{n-1}{2}\rfloor}\vee \lceil \frac{n+1}{2}\rceil{K_1})< \rho(G_{\lceil \frac{2t}{1-t}\rceil})$. 
   	\end{claim}
   	
   	\noindent\emph{Proof of Claim  \ref{claim1}.}~
   	Since $\rho(G_{\lceil \frac{2t}{1-t}\rceil})>\rho(K_{n-\lceil \frac{2t}{1-t}\rceil-1})=n-\lceil \frac{2t}{1-t}\rceil-2$, it suffices to show that $\rho(K_{\lfloor \frac{n-1}{2}\rfloor}\vee \lceil \frac{n+1}{2}\rceil{K_1})\leq n-\lceil \frac{2t}{1-t}\rceil-2$. 
   	Denote $\widetilde{G}=K_{\lfloor \frac{n-1}{2}\rfloor}\vee \lceil \frac{n+1}{2}\rceil{K_1}$.
   	For the partition $V(\widetilde{G})=V(K_{\lfloor \frac{n-1}{2}\rfloor})\bigcup V(\lceil \frac{n+1}{2}\rceil{K_1})$, the quotient matrix of $A(\widetilde{G})$ is
   	\[
   	\widetilde{M}=\begin{pmatrix}
   		\lfloor \frac{n-1}{2}\rfloor-1&\lceil \frac{n+1}{2}\rceil\\
   		\lfloor \frac{n-1}{2}\rfloor&0
   	\end{pmatrix}.
   	\]
   	
   	The characteristic polynomial of $\widetilde{M}$ is
   	\begin{align*}
   		\widetilde{f}(x)=|xI_2-\widetilde{M}|
   		=x(x-\lfloor \frac{n-1}{2}\rfloor+1)-\lfloor \frac{n-1}{2}\rfloor \lceil\frac{n+1}{2}\rceil.
   	\end{align*}
   	
   	We consider the cases of odd and even $n$ separately.
   	When $n$ is odd,
   	\begin{align*}
   		\widetilde{f}(n-\lceil \frac{2t}{1-t}\rceil-2)&=(n-\lceil \frac{2t}{1-t}\rceil-2)(\frac{n}{2}-\lceil \frac{2t}{1-t}\rceil-\frac{1}{2})-\frac{n^2-1}{4}\\
   		&=n(\frac{n}{4}-\frac{3\lceil \frac{2t}{1-t}\rceil+3}{2})+(\lceil \frac{2t}{1-t}\rceil+2)(\lceil \frac{2t}{1-t}\rceil+\frac{1}{2})+\frac{1}{4}.
   	\end{align*}
   	
   	When $n$ is even,
   	\begin{align*}
   		\widetilde{f}(n-\lceil \frac{2t}{1-t}\rceil-2)&=(n-\lceil \frac{2t}{1-t}\rceil-2)(\frac{n}{2}-\lceil \frac{2t}{1-t}\rceil)-\frac{n^2-4}{4}\\
   		&=n(\frac{n}{4}-\frac{3\lceil \frac{2t}{1-t}\rceil+2}{2})+\lceil \frac{2t}{1-t}\rceil(\lceil \frac{2t}{1-t}\rceil+2)+1.
   	\end{align*}
   	
   	Recall that $n\geq \frac{12}{1-t}$. $n\geq 6(\frac{2t}{1-t}+1)+6\geq 6\lceil \frac{2t}{1-t}\rceil+6$. Hence $\widetilde{f}(n-\lceil \frac{2t}{1-t}\rceil-2)>0$,  whether $n$ is odd or even. $\widetilde{f}$ is an upward-opening quadratic function having a negative root. Then $\rho(\widetilde{G})<n-\lceil \frac{2t}{1-t}\rceil-2$. We have $\rho(\widetilde{G})<n-\lceil \frac{2t}{1-t}\rceil-2<\rho(G_{\lceil \frac{2t}{1-t}\rceil})$.
   	This completes the proof of Claim \ref{claim1}. \qed
   	
   	Therefore, in this case, $\rho(G)\leq \rho(K_{\lfloor \frac{n-1}{2}\rfloor}\vee \lceil \frac{n+1}{2}\rceil{K_1}) <\rho(G_{\lceil \frac{2t}{1-t}\rceil})$. This contradicts the assumption $\rho(G)\geq \rho (G_{\lceil \frac{2t}{1-t}\rceil})$.
   	
   	The proof is complete.
    \end{proof}

   \section{Two distance spectral conditions }

	In this section, we establish two tight sufficient conditions in terms of the distance spectral radius for a $t$-tough graph $G$ to have a $\{K_2,\{C_k\}\}$-factor. Theorem  \ref{thm3.5} gives the distance spectral radius condition for $G$ to have a $\{K_2,\{C_k\}\}$-factor, which is divided into three cases based on the value of $t$. Furthermore, for graphs of larger order $n \ge \frac{21-t}{1-t}$, Theorem \ref{thm3} provides a more concise distance spectral radius condition.
	
		\begin{theorem}\label{thm3.5}
		Let $0<t<1$ and let $G$ be a $t$-tough graph of order $n$. Set $G_s=K_s\vee (K_{n-2s-1}+(s+1)K_1)$ with  $\lceil \frac{2t}{1-t}\rceil \leq s \leq \lfloor \frac{n-3}{2}\rfloor$. If one of the following three conditions holds: 
		
		\enspace	(1) $\frac{2\lceil \frac{n}{2}\rceil-2}{2\lfloor  \frac{n}{2}\rfloor+2}< t<1$;
		
		\enspace	(2) $\frac{2\lceil \frac{n}{2}\rceil-4}{2\lceil \frac{n}{2}\rceil}< t\leq \frac{2\lceil \frac{n}{2}\rceil-2}{2\lfloor\frac{n}{2}\rfloor+2}$, $\mu(G)\leq \mu(K_{\lfloor  \frac{n-1}{2}\rfloor}\vee \lceil \frac{n+1}{2}\rceil K_1)$, and $G\neq K_{\lfloor  \frac{n-1}{2}\rfloor}\vee \lceil \frac{n+1}{2}\rceil K_1$;
		
		\enspace	(3) $t\leq \frac{2\lceil \frac{n}{2}\rceil-4}{2\lceil  \frac{n}{2}\rceil}$, $\mu(G)\leq \min\{\mu(K_{\lfloor  \frac{n-1}{2}\rfloor}\vee \lceil \frac{n+1}{2}\rceil K_1), \mu(G_s)\}$, and $G$ is isomorphic to neither $K_{\lfloor  \frac{n-1}{2}\rfloor}\vee \lceil \frac{n+1}{2}\rceil K_1$ nor any $G_s$ with $\left\lceil\frac{2t}{1-t}\right\rceil \leq s\leq \left\lfloor\frac{n-3}{2}\right\rfloor$;
		
		then $G$ has a $\{K_2,\{C_k\}\}$-factor.
	\end{theorem}
	
	\begin{proof}
		The proof of Theorem \ref{thm3.5} is structurally similar to that of Theorem \ref{thm1}; hence, we omit the full proof and place it in Appendix A.
	\end{proof}
	
	\begin{theorem}\label{thm3}
		Let $G$ be a $t$-tough graph of order $n\geq \frac{21-t}{1-t}$, where $0<t<1$. If
		\[
		\mu(G)\leq\mu\left(K_{\lceil\frac{2t}{1-t}\rceil}\vee\left(K_{n-2\lceil\frac{2t}{1-t}\rceil-1}+\left(\lceil\frac{2t}{1-t}\rceil+1\right)K_1\right)\right),
		\]
		then $G$ has a $\{K_2,\{C_k\}\}$-factor unless
		$G= K_{\lceil\frac{2t}{1-t}\rceil}\vee\left(K_{n-2\lceil\frac{2t}{1-t}\rceil-1}+\left(\lceil\frac{2t}{1-t}\rceil+1\right)K_1\right).$
	\end{theorem}
	
	\begin{proof}
		Assume to the contrary that $G$ does not have a $\{K_2,\{C_k\}\}$-factor and is not isomorphic to the exceptional graph stated in the theorem. By Lemma \ref{lem0}, there exists $S \subseteq V(G)$ such that  $i(G-S)-|S|>0$. Let $|S|=s$, $i(G-S)=p$ and $c(G-S)=p+q$. Then $p>s$. We now proceed to discuss the following two cases. 
		
		{\textbf{Case 1.}} $q>0$.
		
		As in Section 3, we can obtain $\mu(G)\geq \mu(G_s)$, where $G_s=K_s\vee(K_{n-2s-1}+(s+1)K_1)$ and $\left\lceil\frac{2t}{1-t}\right\rceil
		\leq s\leq \left\lfloor\frac{n-3}{2}\right\rfloor$. 
		By Lemma \ref{lem8}, we have $\mu(G)\geq \mu(G_s)\geq \mu(G_{\lceil \frac{2t}{1-t}\rceil})$ with equality if and only if $G=G_{\lceil \frac{2t}{1-t}\rceil}$ and $s=\lceil \frac{2t}{1-t}\rceil$. 
		On the other hand, the hypothesis of Theorem \ref{thm3} requires $\mu(G)\leq \mu (G_{\lceil \frac{2t}{1-t}\rceil})$. Consequently,
		$\mu(G)=\mu(G_{\lceil \frac{2t}{1-t}\rceil})$, and hence $G=G_{\lceil \frac{2t}{1-t}\rceil}$, contradicting the assumption that $G$ is not the exceptional graph.
		
		{\textbf{Case 2.}} $q=0$.
		
		In this case, we have $G\preceq K_{n-p}\vee pK_1\preceq K_{\lfloor \frac{n-1}{2}\rfloor}\vee \lceil \frac{n+1}{2}\rceil{K_1}$. Hence $\mu(G)\geq \mu(K_{n-p}\vee p{K_1})\geq \mu(K_{\lfloor \frac{n-1}{2}\rfloor}\vee \lceil \frac{n+1}{2}\rceil{K_1})$. Furthermore, we make the following claim. 
		
		\begin{claim}\label{claim2}
		$\mu(K_{\lfloor\frac{n-1}{2}\rfloor}\vee\lceil \frac{n+1}{2}\rceil K_1)>\mu(G_{\lceil \frac{2t}{1-t}\rceil})$.
		\end{claim}
		
		\noindent\emph{Proof of Claim \ref{claim2}.}~
			Let $a=\lceil\frac{2t}{1-t}\rceil$. Then $n\geq\frac{21-t}{1-t}=10\frac{2t}{1-t}+21>10a+11$, and $n\geq 10a+12$. We first prove that $\mu(K_{\lfloor\frac{n-1}{2}\rfloor}\vee\lceil \frac{n+1}{2}\rceil K_1)>\frac{5n}{4}-1$.
			The equitable quotient matrix of $D(K_{\lfloor\frac{n-1}{2}\rfloor}\vee\lceil \frac{n+1}{2}\rceil K_1)$ is
			\[
			\widetilde{M}
			=
			\begin{pmatrix}
				\lfloor\frac{n-1}{2}\rfloor-1 & \lceil\frac{n+1}{2}\rceil\\
				\lfloor\frac{n-1}{2}\rfloor   & 2\lceil\frac{n+1}{2}\rceil-2
			\end{pmatrix}.
			\]
			Taking $v=(2,3)^{\mathsf T}$ in Lemma \ref{lem11}, we obtain
			\begin{align*}
				\mu(K_{\lfloor\frac{n-1}{2}\rfloor}\vee\lceil \frac{n+1}{2}\rceil K_1)&=\rho(\widetilde{M})\geq\min_i\frac{(\widetilde{M}v)_i}{v_i}\\
				&=\min\left\{\lfloor\frac{n-1}{2}\rfloor-1+\frac{3\lceil\frac{n+1}{2}\rceil}{2},\frac{2\lfloor\frac{n-1}{2}\rfloor}{3}+2\lceil\frac{n+1}{2}\rceil-2\right\}\\
				&\geq\min\left\{\frac{5n-3}{4},\frac{4n-4}{3}\right\}>\frac{5n}{4}-1.
			\end{align*}
			
			Now we  prove that $\mu(G_a)<\frac{5n}{4}-1$. With respect to the partition $V(G_a)=V(K_a)\cup V(K_{n-2a-1})\cup V((a+1)K_1)$, the equitable quotient matrix of $D(G_a)$ is
			\[
			M_a=
			\begin{pmatrix}
				a-1 & n-2a-1   & a+1\\
				a   & n-2a-2 & 2(a+1)\\
				a   & 2(n-2a-1)  & 2a
			\end{pmatrix}.
			\]
			
			Hence $\mu(G_a)=\rho(M_a)$.	Taking $u=(4,4,7)^{\mathsf T}$ in Lemma \ref{lem11}, we obtain
			\[
			\mu(G_a)\leq\max_i\frac{(M_au)_i}{u_i}=\max\left\{n+\frac{3a-1}{4},n+\frac{5a+3}{2},\frac{8n+2a-8}{7}\right\}.
			\]
		    Since $n\geq10a+12$, we have $n+\frac{3a-1}{4}<n+\frac{5a+3}{2}$, $\frac{5n}{4}-1-(n+\frac{5a+3}{2})=\frac{n-10a-10}{4}>0$ and $\frac{5n}{4}-1-\frac{8n+2a-8}{7}=\frac{3n-8a+4}{28}>0$. Consequently, $\mu(G_a)<\frac{5n}{4}-1$. Then $\mu(K_{\lfloor\frac{n-1}{2}\rfloor}\vee\lceil \frac{n+1}{2}\rceil K_1)>\frac{5n}{4}-1>\mu(G_a)$. This completes the proof of Claim \ref{claim2}. \qed

		Therefore, in this case, $\mu(G)\geq\mu\left(K_{\lfloor\frac{n-1}{2}\rfloor}\vee\lceil \frac{n+1}{2}\rceil K_1\right)>\mu(G_{\lceil \frac{2t}{1-t}\rceil})$. This contradicts the assumption $\mu(G)\leq\mu(G_{\lceil \frac{2t}{1-t}\rceil})$.
		
	    The proof is complete.
	\end{proof}

\section{Two distance signless Laplacian spectral conditions }
	
	In this section, we establish two tight sufficient conditions in terms of the distance signless Laplacian spectral radius for a $t$-tough graph $G$ to have a $\{K_2,\{C_k\}\}$-factor. Theorem  \ref{thm4.5} gives the distance signless Laplacian spectral radius condition for $G$ to have a $\{K_2,\{C_k\}\}$-factor, which is divided into three cases based on the value of $t$. Furthermore, for graphs of larger order $n \ge \frac{18+10t}{1-t}$, Theorem \ref{thm4} provides a more concise distance signless Laplacian spectral radius condition.
	
   \begin{theorem}\label{thm4.5}
   	Let $0<t<1$ and let $G$ be a $t$-tough graph of order $n$. Set $G_s=K_s\vee (K_{n-2s-1}+(s+1)K_1)$ with  $\lceil \frac{2t}{1-t}\rceil \leq s \leq \lfloor \frac{n-3}{2}\rfloor$. If one of the following three conditions holds: 
   	
   	\enspace	(1) $\frac{2\lceil \frac{n}{2}\rceil-2}{2\lfloor  \frac{n}{2}\rfloor+2}< t<1$;
   	
   	\enspace	(2) $\frac{2\lceil \frac{n}{2}\rceil-4}{2\lceil \frac{n}{2}\rceil}< t\leq \frac{2\lceil \frac{n}{2}\rceil-2}{2\lfloor\frac{n}{2}\rfloor+2}$, $\eta(G)\leq \eta(K_{\lfloor  \frac{n-1}{2}\rfloor}\vee \lceil \frac{n+1}{2}\rceil K_1)$, and $G\neq K_{\lfloor  \frac{n-1}{2}\rfloor}\vee \lceil \frac{n+1}{2}\rceil K_1$;
   	
   	\enspace	(3) $t\leq \frac{2\lceil \frac{n}{2}\rceil-4}{2\lceil  \frac{n}{2}\rceil}$, $\eta(G)\leq \min\{\eta(K_{\lfloor  \frac{n-1}{2}\rfloor}\vee \lceil \frac{n+1}{2}\rceil K_1), \eta(G_s)\}$, and $G$ is isomorphic to neither $K_{\lfloor  \frac{n-1}{2}\rfloor}\vee \lceil \frac{n+1}{2}\rceil K_1$ nor any $G_s$ with $\left\lceil\frac{2t}{1-t}\right\rceil \leq s\leq \left\lfloor\frac{n-3}{2}\right\rfloor$;
   	
   	then $G$ has a $\{K_2,\{C_k\}\}$-factor.
   \end{theorem}
   
   \begin{proof}
   			The proof of Theorem \ref{thm3.5} is structurally similar to that of Theorem \ref{thm1}; hence, we omit the full proof and place it in Appendix A.
   \end{proof}

	\begin{theorem}\label{thm4}
		Let $G$ be a $t$-tough graph of order $n\geq \frac{18+10t}{1-t}$, where $0<t<1$. If
		\[
		\eta(G)\leq\eta\left(K_{\lceil\frac{2t}{1-t}\rceil}\vee\left(K_{n-2\lceil\frac{2t}{1-t}\rceil-1}+\left(\lceil\frac{2t}{1-t}\rceil+1\right)K_1\right)\right),
		\]
		then $G$ has a $\{K_2,\{C_k\}\}$-factor unless $G= K_{\lceil\frac{2t}{1-t}\rceil}\vee\left(K_{n-2\lceil\frac{2t}{1-t}\rceil-1}+\left(\lceil\frac{2t}{1-t}\rceil+1\right)K_1\right).$
	\end{theorem}
	
	\begin{proof}
		
		Assume to the contrary that $G$ does not have a $\{K_2,\{C_k\}\}$-factor and is not isomorphic to the exceptional graph stated in the theorem. By Lemma \ref{lem0}, there exists $S \subseteq V(G)$ such that  $i(G-S)-|S|>0$. Let $|S|=s$, $i(G-S)=p$ and $c(G-S)=p+q$. Then $p>s$. We now proceed to discuss the following two cases. 
		
		{\textbf{Case 1.}} $q>0$.
		
		As in Section 3, we can obtain $\eta(G)\geq \eta(G_s)$, where $G_s=K_s\vee(K_{n-2s-1}+(s+1)K_1)$ and $\left\lceil\frac{2t}{1-t}\right\rceil\leq s\leq \left\lfloor\frac{n-3}{2}\right\rfloor$. 
		By Lemma \ref{lem9}, we have $\eta(G)\geq \eta(G_s)\geq \eta(G_{\lceil \frac{2t}{1-t}\rceil})$ with equality if and only if $G=G_{\lceil \frac{2t}{1-t}\rceil}$ and $s=\lceil \frac{2t}{1-t}\rceil$. 
		On the other hand, the hypothesis of Theorem \ref{thm4} requires $\eta(G)\leq \eta (G_{\lceil \frac{2t}{1-t}\rceil})$. Consequently, $\eta(G)=\eta(G_{\lceil \frac{2t}{1-t}\rceil})$, and hence $G=G_{\lceil \frac{2t}{1-t}\rceil}$, contradicting the assumption that $G$ is not the exceptional graph.
		
		{\textbf{Case 2.}} $q=0$.
		
		In this case, we have $G\preceq K_{n-p}\vee pK_1\preceq K_{\lfloor \frac{n-1}{2}\rfloor}\vee \lceil \frac{n+1}{2}\rceil{K_1}$. Hence $\eta(G)\geq \eta(K_{n-p}\vee p{K_1})\geq \eta(K_{\lfloor \frac{n-1}{2}\rfloor}\vee \lceil \frac{n+1}{2}\rceil{K_1})$. Furthermore, we make the following claim. 

	    \begin{claim}\label{claim3}
		$\eta(K_{\lfloor\frac{n-1}{2}\rfloor}\vee\lceil \frac{n+1}{2}\rceil K_1)>\eta(G_{\lceil \frac{2t}{1-t}\rceil})$.
		\end{claim}
		
		\noindent\emph{Proof of Claim \ref{claim3}.}~
		Let $a=\lceil \frac{2t}{1-t}\rceil$. Then $n\geq\frac{18+10t}{1-t}=14\frac{2t}{1-t}+18>14a+4$, and $n\geq 14a+5$.
		We first prove that $\eta(K_{\lfloor\frac{n-1}{2}\rfloor}\vee\lceil\frac{n+1}{2}\rceil K_1)>\frac{8(n-1)}{3}$. The equitable quotient matrix of $Q(K_{\lfloor\frac{n-1}{2}\rfloor}\vee\lceil\frac{n+1}{2}\rceil K_1)$ is
	    \[
	    \widetilde{N}=\begin{pmatrix} n+\lfloor\frac{n-1}{2}\rfloor-2 & \lceil\frac{n+1}{2}\rceil\\ \lfloor\frac{n-1}{2}\rfloor & n+3\lceil\frac{n+1}{2}\rceil-4 \end{pmatrix}.
	    \]
		Hence $\eta(K_{\lfloor\frac{n-1}{2}\rfloor}\vee\lceil\frac{n+1}{2}\rceil K_1)=\rho(\widetilde{N})$. Taking $v=(1,3)^{\mathsf T}$, we obtain $\frac{(\widetilde{N}v)_1}{v_1}=4n-2\lfloor\frac{n-1}{2}\rfloor-2\geq3n-1>\frac{8(n-1)}{3}$ and $\frac{(\widetilde{N}v)_2}{v_2}=4n-\frac{8\lfloor\frac{n-1}{2}\rfloor}{3}-4\geq\frac{8(n-1)}{3}$. Thus $\widetilde{N}v\geq\frac{8(n-1)}{3}v$, with a strict inequality in the first coordinate. Since $\widetilde{N}$ is irreducible and $v>0$, Lemma \ref{lem11} yields $\eta(K_{\lfloor\frac{n-1}{2}\rfloor}\vee\lceil\frac{n+1}{2}\rceil K_1)>\frac{8(n-1)}{3}$.
			
		Now we prove that $\eta(G_a)<\frac{8(n-1)}{3}$. With respect to the partition $V(G_a)=V(K_a)\cup V(K_{n-2a-1})\cup V((a+1)K_1)$, the equitable quotient matrix of $Q(G_a)$ is
		\[N_a=\begin{pmatrix} n+a-2 & n-2a-1 & a+1\\ a & 2n-a-2 & 2(a+1)\\ a & 2(n-2a-1) & 2n+a-2 \end{pmatrix}.\]
		Hence $\eta(G_a)=\rho(N_a)$. Taking $u=(1,1,3)^{\mathsf T}$, we obtain $\frac{(N_au)_1}{u_1}=2n+2a$, $\frac{(N_au)_2}{u_2}=2n+6a+4$ and $\frac{(N_au)_3}{u_3}=\frac{8(n-1)}{3}$. Since $n\geq9a+10$, we have $\frac{8(n-1)}{3}-(2n+2a)=\frac{2n-6a-8}{3}>0$ and $\frac{8(n-1)}{3}-(2n+6a+4)=\frac{2n-18a-20}{3}\geq0$. Therefore, $N_au\leq\frac{8(n-1)}{3}u$, with a strict inequality in the first coordinate. Since $N_a$ is irreducible and $u>0$, Lemma \ref{lem11} implies $\eta(G_a)<\frac{8(n-1)}{3}$.
			
		In this case, $\eta(K_{\lfloor\frac{n-1}{2}\rfloor}\vee\lceil\frac{n+1}{2}\rceil K_1)>\frac{8(n-1)}{3}>\eta(G_a)$. It remains to consider $n<9a+10$. Since $n\geq14a+4$ and $a\geq1$ is an integer, we obtain that $a=1$ and $n=18$. For the same vectors $u$ and $v$, Lemma \ref{lem11} gives $\eta(G_a)\leq\max\{38,46,\frac{136}{3}\}=46$ and $\eta(K_{\lfloor\frac{n-1}{2}\rfloor}\vee\lceil\frac{n+1}{2}\rceil K_1)\geq\min\{54,\frac{140}{3}\}=\frac{140}{3}>46$. Thus the desired inequality also holds in this case. This completes the proof of Claim \ref{claim3}. \qed

		Therefore, in this case, $\eta(G)\geq\eta\left(K_{\lfloor\frac{n-1}{2}\rfloor}\vee\lceil \frac{n+1}{2}\rceil K_1\right)>\eta(G_{\lceil \frac{2t}{1-t}\rceil})$. This contradicts the assumption $\eta(G)\leq\eta(G_{\lceil \frac{2t}{1-t}\rceil})$.
		
	    The proof is complete.
	\end{proof}

	\section*{Declaration of competing interest}\label{A}
	
	The authors declare that they have no known competing financial interests or personal relationships that could have appeared to influence the work reported in this paper.
	
	\section*{Data availability}
	
	No data was used for the research described in the article.

	\section*{Appendix A: Proofs of Theorems~\ref{thm3.5} and~\ref{thm4.5}}
	\label{appen}

	\textbf{Proof of Theorem \ref{thm3.5}}
	
	Assume to the contrary that $G$ does not have a $\{K_2,\{C_k\}\}$-factor. By Lemma \ref{lem0}, there exists $S \subseteq V(G)$ such that  $i(G-S)-|S|>0$. Let $|S|=s$, $i(G-S)=p$ and $c(G-S)=p+q$. Then $p>s$. We now proceed to discuss the following three cases. 
	
	{\textbf{Case 1.}} $\frac{2\lceil \frac{n}{2}\rceil-2}{2\lfloor  \frac{n}{2}\rfloor+2}< t<1$.
	
	If there are nontrivial components in $G-S$, we can add all possible edges between the nontrivial components of $G-S$ so that they form a clique. Then $G\preceq K_s\vee (K_{n-s-p}+p{K_1})$. Set $H_1=K_s\vee (K_{n-s-p}+p{K_1})$. Then $t(H_1)\geq t(G)\geq t$ and $H_1$ does not have a $\{K_2,\{C_k\}\}$-factor. By Lemma \ref{lem10}(1), $t\leq \frac{2\lceil \frac{n}{2}\rceil-4}{2\lceil  \frac{n}{2}\rceil}$. This implies that $t\leq \frac{2\lceil \frac{n}{2}\rceil-4}{2\lceil  \frac{n}{2}\rceil}< \frac{2\lceil \frac{n}{2}\rceil-2}{2\lfloor  \frac{n}{2}\rfloor+2}< t$, which is a contradiction.
	
	If there are no nontrivial components in $G-S$, then $G \preceq K_{n-p}\vee pK_1$. Set $H_2=K_{n-p}\vee pK_1$. Since $t(H_2)\geq t(G)\geq t$ and $H_2$ does not have a $\{K_2,\{C_k\}\}$-factor, by Lemma \ref{lem10}(2), $t\leq \frac{2\lceil \frac{n}{2}\rceil-2}{2\lfloor\frac{n}{2}\rfloor+2}$. This implies that $t\leq \frac{2\lceil \frac{n}{2}\rceil-2}{2\lfloor\frac{n}{2}\rfloor+2}<t$, which is a contradiction. 
	
	Thus, both possibilities lead to a contradiction. The proof of Theorem \ref{thm3.5}(1) is complete.

	{\textbf{Case 2.}} $\frac{2\lceil \frac{n}{2}\rceil-4}{2\lceil \frac{n}{2}\rceil}< t\leq \frac{2\lceil \frac{n}{2}\rceil-2}{2\lfloor\frac{n}{2}\rfloor+2}$.
	
	If there are nontrivial components in $G-S$, $G\preceq K_s\vee (K_{n-s-p}+p{K_1})$. Set $H_1=K_s\vee (K_{n-s-p}+p{K_1})$. Then $t(H_1)\geq t(G)\geq t$ and $H_1$ does not have a $\{K_2,\{C_k\}\}$-factor. By Lemma \ref{lem10}(1), $t\leq \frac{2\lceil \frac{n}{2}\rceil-4}{2\lceil  \frac{n}{2}\rceil}$. This implies that $t\leq \frac{2\lceil \frac{n}{2}\rceil-4}{2\lceil  \frac{n}{2}\rceil}<t$, which is a contradiction.
	
	If there are no nontrivial components in $G-S$, then $G \preceq K_{n-p}\vee pK_1$. Since $p>s=n-p$, we have $p\geq \lceil \frac{n+1}{2}\rceil$. By moving $p-\lceil \frac{n+1}{2}\rceil$ vertices from the independent set into the clique and adding the required edges, we see that $K_{n-p}\vee pK_1\preceq K_{\lfloor  \frac{n-1}{2}\rfloor}\vee \lceil \frac{n+1}{2}\rceil K_1$. Therefore, $\mu(G)\geq \mu(K_{n-p}\vee pK_1)\geq \mu(K_{\lfloor \frac{n-1}{2}\rfloor}\vee \lceil \frac{n+1}{2}\rceil K_1)$ with equality throughout if and only if $G=K_{\lfloor  \frac{n-1}{2}\rfloor}\vee \lceil \frac{n+1}{2}\rceil K_1$. Together with the hypothesis $\mu(G)\leq\mu(K_{\lfloor  \frac{n-1}{2}\rfloor}\vee \lceil \frac{n+1}{2}\rceil K_1)$, this forces $G=K_{\lfloor  \frac{n-1}{2}\rfloor}\vee \lceil \frac{n+1}{2}\rceil K_1$, contradicting the exclusion in (2). 
	
	The proof of Theorem \ref{thm3.5}(2) is complete.
	
	{\textbf{Case 3.}} $t\leq \frac{2\lceil \frac{n}{2}\rceil-4}{2\lceil  \frac{n}{2}\rceil}$.
	
	If there are nontrivial components in $G-S$, $G\preceq K_s\vee (K_{n-s-p}+p{K_1})$ and $n-s-p\geq 2$. Since $G$ is $t$-tough and $c(G-S)=p+q\geq p+1$,
	we have $s\geq t(p+q)\geq t(p+1)\geq t(s+2)$. Together with $p\geq s+1$ and $n-s-p\geq 2$, it gives that $\lceil\frac{2t}{1-t}\rceil\leq s\leq\lfloor\frac{n-3}{2}\rfloor$. By moving $p-s-1$ isolated vertices into the clique $K_{n-s-p}$ and adding the required edges, we obtain $H_1\preceq G_s=K_s\vee (K_{n-2s-1}+(s+1)K_1)$. Hence $\mu(G)\geq \mu(H_1)\geq \mu(G_s)\geq \min\{\mu(K_{\lfloor  \frac{n-1}{2}\rfloor}\vee \lceil \frac{n+1}{2}\rceil K_1), \mu(G_s)\}$. If $\mu(G)>\mu(K_{\lfloor  \frac{n-1}{2}\rfloor}\vee \lceil \frac{n+1}{2}\rceil K_1)$, then $\mu(G)>\min\{\mu(K_{\lfloor  \frac{n-1}{2}\rfloor}\vee \lceil \frac{n+1}{2}\rceil K_1), \mu(G_s)\}$, contradicting the hypothesis. If $\mu(G_s)=\min\{\mu(K_{\lfloor  \frac{n-1}{2}\rfloor}\vee \lceil \frac{n+1}{2}\rceil K_1), \mu(G_s)\}$, the hypothesis $\mu(G)\leq \min\{\mu(K_{\lfloor  \frac{n-1}{2}\rfloor}\vee \lceil \frac{n+1}{2}\rceil K_1), \mu(G_s)\}$ forces equality throughout. Then $G= G_s$ or $G=K_{\lfloor  \frac{n-1}{2}\rfloor}\vee \lceil \frac{n+1}{2}\rceil K_1$, contrary to the exclusion in (3).
	Thus, both possibilities lead to a contradiction.
	
	If there are no nontrivial components in $G-S$, then $G \preceq K_{n-p}\vee pK_1\preceq K_{\lfloor  \frac{n-1}{2}\rfloor}\vee \lceil \frac{n+1}{2}\rceil K_1$. Hence $\mu(G)\geq \mu(K_{\lfloor  \frac{n-1}{2}\rfloor}\vee \lceil \frac{n+1}{2}\rceil K_1)\geq \min\{\mu(K_{\lfloor  \frac{n-1}{2}\rfloor}\vee \lceil \frac{n+1}{2}\rceil K_1),  \mu(G_s)\}$. If $\mu(G)>\mu(K_{\lfloor  \frac{n-1}{2}\rfloor}\vee \lceil \frac{n+1}{2}\rceil K_1)$, then $\mu(G)>\min\{\mu(K_{\lfloor  \frac{n-1}{2}\rfloor}\vee \lceil \frac{n+1}{2}\rceil K_1), \mu(G_s)\}$, contradicting the hypothesis. If $\mu(G)=\min\{\mu(K_{\lfloor  \frac{n-1}{2}\rfloor}\vee \lceil \frac{n+1}{2}\rceil K_1), \mu(G_s)\}$, the hypothesis $\mu(G)\leq \min\{\mu(K_{\lfloor  \frac{n-1}{2}\rfloor}\vee \lceil \frac{n+1}{2}\rceil K_1), \mu(G_s)\}$ forces equality throughout. Then $G= G_s$ or $G=K_{\lfloor  \frac{n-1}{2}\rfloor}\vee \lceil \frac{n+1}{2}\rceil K_1$, contrary to the exclusion in (3).
	Thus, both possibilities lead to a contradiction.
	
	The proof of Theorem \ref{thm3.5}(3) is complete.	\qed
	
	\textbf{Proof of Theorem \ref{thm4.5}}
	
	Assume to the contrary that $G$ does not have a $\{K_2,\{C_k\}\}$-factor. By Lemma \ref{lem0}, there exists $S \subseteq V(G)$ such that  $i(G-S)-|S|>0$. Let $|S|=s$, $i(G-S)=p$ and $c(G-S)=p+q$. Then $p>s$. We now proceed to discuss the following three cases. 
	
	{\textbf{Case 1.}} $\frac{2\lceil \frac{n}{2}\rceil-2}{2\lfloor  \frac{n}{2}\rfloor+2}< t<1$.
	
	If there are nontrivial components in $G-S$, we can add all possible edges between the nontrivial components of $G-S$ so that they form a clique. Then $G\preceq K_s\vee (K_{n-s-p}+p{K_1})$. Set $H_1=K_s\vee (K_{n-s-p}+p{K_1})$. Then $t(H_1)\geq t(G)\geq t$ and $H_1$ does not have a $\{K_2,\{C_k\}\}$-factor. By Lemma \ref{lem10}(1), $t\leq \frac{2\lceil \frac{n}{2}\rceil-4}{2\lceil  \frac{n}{2}\rceil}$. This implies that $t\leq \frac{2\lceil \frac{n}{2}\rceil-4}{2\lceil  \frac{n}{2}\rceil}< \frac{2\lceil \frac{n}{2}\rceil-2}{2\lfloor  \frac{n}{2}\rfloor+2}< t$, which is a contradiction.
	
	If there are no nontrivial components in $G-S$, then $G \preceq K_{n-p}\vee pK_1$. Set $H_2=K_{n-p}\vee pK_1$. Since $t(H_2)\geq t(G)\geq t$ and $H_2$ does not have a $\{K_2,\{C_k\}\}$-factor, by Lemma \ref{lem10}(2), $t\leq \frac{2\lceil \frac{n}{2}\rceil-2}{2\lfloor\frac{n}{2}\rfloor+2}$. This implies that $t\leq \frac{2\lceil \frac{n}{2}\rceil-2}{2\lfloor\frac{n}{2}\rfloor+2}<t$, which is a contradiction. 
	
	Thus, both possibilities lead to a contradiction. The proof of Theorem \ref{thm4.5}(1) is complete.

	{\textbf{Case 2.}} $\frac{2\lceil \frac{n}{2}\rceil-4}{2\lceil \frac{n}{2}\rceil}< t\leq \frac{2\lceil \frac{n}{2}\rceil-2}{2\lfloor\frac{n}{2}\rfloor+2}$.
	
	If there are nontrivial components in $G-S$, $G\preceq K_s\vee (K_{n-s-p}+p{K_1})$. Set $H_1=K_s\vee (K_{n-s-p}+p{K_1})$. Then $t(H_1)\geq t(G)\geq t$ and $H_1$ does not have a $\{K_2,\{C_k\}\}$-factor. By Lemma \ref{lem10}(1), $t\leq \frac{2\lceil \frac{n}{2}\rceil-4}{2\lceil  \frac{n}{2}\rceil}$. This implies that $t\leq \frac{2\lceil \frac{n}{2}\rceil-4}{2\lceil  \frac{n}{2}\rceil}<t$, which is a contradiction.
	
	If there are no nontrivial components in $G-S$, then $G \preceq K_{n-p}\vee pK_1$. Since $p>s=n-p$, we have $p\geq \lceil \frac{n+1}{2}\rceil$. By moving $p-\lceil \frac{n+1}{2}\rceil$ vertices from the independent set into the clique and adding the required edges, we see that $K_{n-p}\vee pK_1\preceq K_{\lfloor  \frac{n-1}{2}\rfloor}\vee \lceil \frac{n+1}{2}\rceil K_1$. Therefore, $\eta(G)\geq \eta(K_{n-p}\vee pK_1)\geq \eta(K_{\lfloor \frac{n-1}{2}\rfloor}\vee \lceil \frac{n+1}{2}\rceil K_1)$ with equality throughout if and only if $G=K_{\lfloor  \frac{n-1}{2}\rfloor}\vee \lceil \frac{n+1}{2}\rceil K_1$. Together with the hypothesis $\eta(G)\leq\eta(K_{\lfloor  \frac{n-1}{2}\rfloor}\vee \lceil \frac{n+1}{2}\rceil K_1)$, this forces $G=K_{\lfloor  \frac{n-1}{2}\rfloor}\vee \lceil \frac{n+1}{2}\rceil K_1$, contradicting the exclusion in (2). 
	
	The proof of Theorem \ref{thm4.5}(2) is complete.
	
	{\textbf{Case 3.}} $t\leq \frac{2\lceil \frac{n}{2}\rceil-4}{2\lceil  \frac{n}{2}\rceil}$.
	
	If there are nontrivial components in $G-S$, $G\preceq K_s\vee (K_{n-s-p}+p{K_1})$ and $n-s-p\geq 2$. Since $G$ is $t$-tough and $c(G-S)=p+q\geq p+1$,
	we have $s\geq t(p+q)\geq t(p+1)\geq t(s+2)$. Together with $p\geq s+1$ and $n-s-p\geq 2$, it gives that $\lceil\frac{2t}{1-t}\rceil\leq s\leq\lfloor\frac{n-3}{2}\rfloor$. By moving $p-s-1$ isolated vertices into the clique $K_{n-s-p}$ and adding the required edges, we obtain $H_1\preceq G_s=K_s\vee (K_{n-2s-1}+(s+1)K_1)$. Hence $\eta(G)\geq \eta(H_1)\geq \eta(G_s)\geq \min\{\eta(K_{\lfloor  \frac{n-1}{2}\rfloor}\vee \lceil \frac{n+1}{2}\rceil K_1), \eta(G_s)\}$. If $\eta(G)>\eta(K_{\lfloor  \frac{n-1}{2}\rfloor}\vee \lceil \frac{n+1}{2}\rceil K_1)$, then $\eta(G)>\min\{\eta(K_{\lfloor  \frac{n-1}{2}\rfloor}\vee \lceil \frac{n+1}{2}\rceil K_1), \eta(G_s)\}$, contradicting the hypothesis. If $\eta(G_s)=\min\{\eta(K_{\lfloor  \frac{n-1}{2}\rfloor}\vee \lceil \frac{n+1}{2}\rceil K_1), \eta(G_s)\}$, the hypothesis $\eta(G)\leq \min\{\eta(K_{\lfloor  \frac{n-1}{2}\rfloor}\vee \lceil \frac{n+1}{2}\rceil K_1), \eta(G_s)\}$ forces equality throughout. Then $G= G_s$ or $G=K_{\lfloor  \frac{n-1}{2}\rfloor}\vee \lceil \frac{n+1}{2}\rceil K_1$, contrary to the exclusion in (3).
	Thus, both possibilities lead to a contradiction.
	
	If there are no nontrivial components in $G-S$, then $G \preceq K_{n-p}\vee pK_1\preceq K_{\lfloor  \frac{n-1}{2}\rfloor}\vee \lceil \frac{n+1}{2}\rceil K_1$. Hence $\eta(G)\geq \eta(K_{\lfloor  \frac{n-1}{2}\rfloor}\vee \lceil \frac{n+1}{2}\rceil K_1)\geq \min\{\eta(K_{\lfloor  \frac{n-1}{2}\rfloor}\vee \lceil \frac{n+1}{2}\rceil K_1),  \eta(G_s)\}$. If $\eta(G)>\eta(K_{\lfloor  \frac{n-1}{2}\rfloor}\vee \lceil \frac{n+1}{2}\rceil K_1)$, then $\eta(G)>\min\{\eta(K_{\lfloor  \frac{n-1}{2}\rfloor}\vee \lceil \frac{n+1}{2}\rceil K_1), \eta(G_s)\}$, contradicting the hypothesis. If $\eta(G)=\min\{\eta(K_{\lfloor  \frac{n-1}{2}\rfloor}\vee \lceil \frac{n+1}{2}\rceil K_1), \eta(G_s)\}$, the hypothesis $\eta(G)\leq \min\{\eta(K_{\lfloor  \frac{n-1}{2}\rfloor}\vee \lceil \frac{n+1}{2}\rceil K_1), \eta(G_s)\}$ forces equality throughout. Then $G= G_s$ or $G=K_{\lfloor  \frac{n-1}{2}\rfloor}\vee \lceil \frac{n+1}{2}\rceil K_1$, contrary to the exclusion in (3).
	Thus, both possibilities lead to a contradiction.
	
	The proof of Theorem \ref{thm4.5}(3) is complete.	\qed

\end{document}